\documentclass[a4paper]{amsart}

\usepackage{amsmath}
\usepackage{amssymb}
\usepackage{amsthm}
\usepackage{enumitem}
\usepackage{mathtools}
\usepackage{tikz-cd}
\usepackage{xcolor}

\usepackage{hyperref}

\newtheorem*{thm*}{Theorem}
\newtheorem{thm}{Theorem}[section]
\newtheorem{lem}[thm]{Lemma}

\theoremstyle{definition}
\newtheorem{df}[thm]{Definition}
\newtheorem{ex}[thm]{Example}
\theoremstyle{remark}
\newtheorem{rmk}[thm]{Remark}

\numberwithin{equation}{section}

\makeatletter
\g@addto@macro\bfseries{\boldmath}
\makeatother

\DeclareMathOperator{\Aut}{Aut}   
\DeclareMathOperator{\Coin}{Coin} 
\DeclareMathOperator{\coin}{coin} 
\DeclareMathOperator{\Fix}{Fix}   
\newcommand{\g}{\mathfrak{g}}     
\DeclareMathOperator{\GL}{GL}	  
\DeclareMathOperator{\id}{id}	  
\DeclareMathOperator{\ind}{ind}   
\newcommand{\orb}{\backslash}     
\newcommand{\R}{\mathcal{R}}      
\DeclareMathOperator{\RT}{RT}     
\renewcommand{\S}{\Sigma}         
\DeclareMathOperator{\tr}{tr}	  

\newcommand{\FF}{\mathbb{F}}      
\newcommand{\QQ}{\mathbb{Q}}      
\newcommand{\RR}{\mathbb{R}}      
\newcommand{\ZZ}{\mathbb{Z}}      

\newcommand{\eps}{\varepsilon}
\renewcommand{\iff}{\enskip\Leftrightarrow\enskip}

\title[Averaging formulas for $n$-valued maps]{Averaging formulas for the Reidemeister trace, Lefschetz and Nielsen numbers of $n$-valued maps}
\author{Karel Dekimpe and Lore De Weerdt}
\thanks{Research supported by Methusalem grant METH/21/03 -- long term structural funding of the Flemish Government.}
\address{KU Leuven Campus Kulak Kortrijk, 8500 Kortrijk, Belgium}
\email{Karel.Dekimpe@kuleuven.be}
\email{Lore.DeWeerdt@kuleuven.be}

\begin{document}

\maketitle

\vspace*{-5mm}

\begin{abstract}
For an $n$-valued self-map $f$ of a closed manifold $X$, we prove an averaging formula for the Reidemeister trace of $f$ in terms of the Reidemeister coincidence traces of single-valued maps between finite orientable covering spaces of $X$. We then derive analogous formulas for the Lefschetz and Nielsen numbers of $f$. In the special case where $X$ is an infra-nilmanifold, we obtain explicit formulas for the Lefschetz and Nielsen numbers of any $n$-valued map on $X$.
\end{abstract}

\section{Introduction}

If $X$ is a compact polyhedron, the Lefschetz number of a continuous self-map $f:X\to X$ is \begin{equation}\label{eq:L(f)}
L(f)=\sum_k (-1)^k \tr(f_*:H_k(X,\QQ)\to H_k(X,\QQ)).
\end{equation}
The Lefschetz fixed point theorem (see e.g.\ \cite{hatcher}) states that, if $L(f)\neq 0$, then any map homotopic to $f$ has at least one fixed point.

More precise information about the existence of fixed points is given by the Nielsen number $N(f)$. This number is defined by subdividing the fixed point set $\Fix(f)=\{x\in X\mid f(x)=x\}$ into so-called fixed point classes, and assigning an index to each fixed point class. The Nielsen number is the number of fixed point classes with non-zero index. It is a lower bound for the number of fixed points among all maps homotopic to $f$ (which is sharp when $X$ is a manifold of dimension at least $3$; see \cite{wecken3}). By the Lefschetz-Hopf formula, the Lefschetz number can equivalently be computed as the sum of the indices of all fixed point classes.

Algebraic analogues of the fixed point classes are the Reidemeister classes (or twisted conjugacy classes) of $\pi_1(X)$ which are determined by the induced morphism of $f$ on $\pi_1(X)$. The Reidemeister trace of $f$, originally defined by Reidemeister in \cite{reidemeister1936}, is a certain element of the free abelian group generated by these Reidemeister classes, i.e.\ an integral linear combination of Reidemeister classes.
In \cite{wecken2}, Wecken showed that the coefficient of each Reidemeister class in this linear combination is precisely the index of the corresponding fixed point class.

\medskip

If $M_1$ and $M_2$ are closed orientable manifolds of equal dimension, also the Lefschetz and Nielsen coincidence numbers and the coincidence Reidemeister trace of a pair of maps $f,g:M_1\to M_2$ can be defined. These are analogous invariants containing information about the coincidence set $\Coin(f,g)=\{x\in M_1 \mid f(x)=g(x) \}$.

Many averaging formulas for these invariants in terms of lifts to a finite covering space have appeared in the literature, both in fixed point theory and in coincidence theory. If $\bar{X}$ is a finite cover of $X$ with covering group $F$, and $\bar{f}:\bar{X}\to \bar{X}$ is a lift of $f$, then the set of all lifts of $f$ to $\bar{X}$ is $\{\bar{\alpha}\bar{f}\mid \bar{\alpha}\in F \}$, and \[
L(f)=\frac{1}{|F|}\sum_{\bar{\alpha}\in F} L(\bar{\alpha}\bar{f})
\]
(see e.g.\ \cite{jiang}). Similar formulas hold for the Lefschetz coincidence number and the fixed point and coincidence Reidemeister traces (see \cite{leestaecker}). For the Nielsen number, in general only the inequality \[
N(f)\geq \frac{1}{|F|}\sum_{\bar{\alpha}\in F} N(\bar{\alpha}\bar{f})
\] 
can be obtained. A necessary and sufficient condition for the equality is given in \cite{kimleelee}, and this condition is shown to hold e.g.\ when $X$ is an infra-nilmanifold and $\bar{X}$ is a nilmanifold. Similar results for coincidences are shown in \cite{kimlee2007}.

\medskip

Recently, the Lefschetz and Nielsen numbers and the Reidemeister trace have been extended to fixed point theory of $n$-valued maps (see \cite{brown2007,schirmer2,staecker2023}). A straightforward generalisation of the single-valued averaging formulas to this setting would be possible, but would not be very useful. For instance in the setting of infra-nilmanifolds, single-valued fixed point theory is very well understood: any single-valued map on an infra-nilmanifold can be lifted to a finite cover which is a nilmanifold, and the Lefschetz and Nielsen numbers of these lifts can be computed using a homotopy with a linear map. In the $n$-valued case, both steps of this process fail: $n$-valued maps in general do not admit lifts to a nilmanifold, so no averaging formula can be applied; and Lefschetz and Nielsen numbers of $n$-valued maps on nilmanifolds can no longer be computed using a homotopy with a linear map, since $n$-valued maps are in general not homotopic to linear-like maps (see \cite{non-affine}).

In this paper, we link fixed points of $n$-valued maps to coincidences of single-valued maps, where the single-valued maps go from one finite covering space of $X$ to another. We prove an averaging formula for the Reidemeister trace of an $n$-valued map in terms of the coincidence Reidemeister traces of these associated single-valued maps, and we derive similar averaging formulas for the Lefschetz and Nielsen numbers. In the special case of infra-nilmanifolds, we can plug in the known formulas for the single-valued Lefschetz and Nielsen coincidence numbers on nilmanifolds, to obtain practical formulas for the Lefschetz and Nielsen numbers of an $n$-valued map.

\section{Lefschetz number, Nielsen number and Reidemeister trace}

Let $X$ be a compact polyhedron with universal cover $p:\tilde{X}\to X$, and let $\pi$ be the group of covering transformations of $p$ (which is isomorphic to the fundamental group of $X$). A map $f:X\to X$ with lift $\tilde{f}:\tilde{X}\to\tilde{X}$ induces a morphism $\varphi:\pi\to\pi$ defined by \[
\tilde{f}\alpha = \varphi(\alpha)\tilde{f}, \quad \forall\alpha\in \pi.
\]
The set of Reidemeister classes $\R[\varphi]$ is the quotient of $\pi$ by the (twisted conjugacy) relation \[
\alpha \sim_\varphi \beta \iff \exists \gamma\in \pi: \alpha=\gamma\beta\varphi(\gamma^{-1}).
\]
We will write the Reidemeister class of $\alpha$ as $[\alpha]_\varphi$, or $[\alpha]$ when the context is clear.

\begin{rmk}
Any other lift of $f$ can be written as $\alpha\tilde{f}$ for some $\alpha\in \pi$. The morphism induced by the lift $\alpha\tilde{f}$ is $\tau_\alpha\varphi:\pi\to\pi$, where $\tau_\alpha$ is conjugation with $\alpha$.
\end{rmk}

The Reidemeister classes stand in one-to-one correspondence with the Nielsen fixed point classes $p\Fix(\alpha\tilde{f})$. For each Nielsen class, let $\ind(f;p\Fix(\alpha\tilde{f}))$ denote its fixed point index. Then the Lefschetz number and the Nielsen number of $f$ can be written as \begin{align*}
L(f)&=\sum_{[\alpha]\in \R[\varphi]} \ind(f;p\Fix(\alpha\tilde{f})) \\
N(f)&=\#\{[\alpha]\in \R[\varphi] \mid \ind(f;p\Fix(\alpha\tilde{f}))\neq 0 \},
\end{align*}
and the Reidemeister trace as \[
\RT(f,\tilde{f})=\sum_{[\alpha]\in \R[\varphi]} \ind(f;p\Fix(\alpha\tilde{f}))\,[\alpha]
\]
in the free abelian group $\ZZ\R[\varphi]$ generated by the Reidemeister classes. Like the trace formula \eqref{eq:L(f)} for the Lefschetz number, there also exists a trace formula for the Reidemeister trace; see e.g.\ Reidemeister's original definition in \cite{reidemeister1936}.

Note that the Reidemeister trace combines the information of the Lefschetz and Nielsen numbers: the sum of the coefficients of $\RT(f,\tilde{f})$ is equal to the Lefschetz number of $f$, and the number of non-zero terms of $\RT(f,\tilde{f})$ is equal to the Nielsen number. 

Suppose $\Gamma$ is a finite index normal subgroup of $\pi$ so that $\varphi(\Gamma)\subseteq \Gamma$, and let $\bar{X}=\Gamma\orb \tilde{X}$ be the corresponding finite cover of $X$. Then there exists a lift $\bar{f}:\bar{X}\to \bar{X}$ of $f$ so that $\tilde{f}:\tilde{X}\to\tilde{X}$ is also a lift of $\bar{f}$, i.e.\ so that the following diagram commutes: \[
\begin{tikzcd}
\tilde{X} \ar[d] \ar[r,"\tilde{f}"] & \tilde{X} \ar[d] \\
\bar{X} \ar[d] \ar[r,"\bar{f}"] & \bar{X} \ar[d] \\
X \ar[r,"\tilde{f}"] & X.
\end{tikzcd}
\] 
For any $\alpha\in \pi$, the map $\alpha\tilde{f}$ is a lift of $\bar{\alpha}\bar{f}$, where $\bar{\alpha}\in \pi/\Gamma$ is the projection of $\alpha$. If $\varphi':\Gamma\to \Gamma$ denotes the morphism induced by $\varphi$, the corresponding morphism of covering groups induced by $\bar{\alpha}\bar{f}$ is $\tau_\alpha\varphi'$. There is a natural map $\hat{r}^\alpha:\ZZ\R[\tau_\alpha\varphi']\to \ZZ\R[\varphi]$ given by $\hat{r}^\alpha([\beta]_{\tau_\alpha\varphi'})=[\beta\alpha]_\varphi$. With these notations, we can state the averaging formula for the Reidemeister trace:

\begin{thm}[{\cite[Theorem 1.3]{leestaecker}}]
For a map $f:X\to X$ with lift $\bar{f}:\bar{X}\to \bar{X}$ as above, \[
\RT(f,\tilde{f})=\frac{1}{[\pi:\Gamma]} \sum_{\bar{\alpha}\in \pi/\Gamma} \hat{r}^\alpha(\RT(\bar{\alpha}\bar{f},\alpha\tilde{f})).
\]
\end{thm}

The Lefschetz and Nielsen numbers and the Reidemeister trace can be generalised to coincidence theory as follows. Let $M_1$ and $M_2$ be equal-dimensional closed orientable manifolds with universal covers $p_1:\tilde{M}_1\to M_1$ and $p_2:\tilde{M}_2\to M_2$, with covering groups $\pi_1$ and $\pi_2$. For two maps $f,g:M_1\to M_2$ with lifts $\tilde{f},\tilde{g}:\tilde{M}_1\to\tilde{M}_2$, let $\varphi,\psi:\pi_1\to \pi_2$ be the induced morphisms defined by \[
\tilde{f}\alpha = \varphi(\alpha)\tilde{f},\enskip\tilde{g}\alpha = \psi(\alpha)\tilde{g}, \quad \forall\alpha\in \pi_1.
\]
The set of Reidemeister classes $\R[\varphi,\psi]$\label{R-coin} is the quotient of $\pi_2$ by the (doubly twisted conjugacy) relation \[
\alpha \sim_{\varphi,\psi} \beta \iff \exists \gamma\in \pi_1: \alpha=\psi(\gamma)\beta\varphi(\gamma^{-1}).
\]
We will write the Reidemeister class of $\alpha$ as $[\alpha]_{\varphi,\psi}$, or $[\alpha]$ when the context is clear.

The Reidemeister classes partition the coincidence set of $f$ and $g$ into coincidence classes $p_1\!\Coin(\alpha\tilde{f},\tilde{g})$. Let $\ind(f,g;p_1\!\Coin(\alpha\tilde{f},\tilde{g}))$ denote the coincidence index of a coincidence class. Then the Lefschetz and Nielsen coincidence numbers and the coincidence Reidemeister trace are \begin{align*}
L(f,g)&=\sum_{[\alpha]\in \R[\varphi,\psi]} \ind(f,g;p_1\!\Coin(\alpha\tilde{f},\tilde{g})) \\
N(f,g)&=\#\{[\alpha]\in \R[\varphi,\psi] \mid \ind(f,g;p_1\!\Coin(\alpha\tilde{f},\tilde{g}))\neq 0 \} \\
\RT(f,\tilde{f};g,\tilde{g})&=\sum_{[\alpha]\in \R[\varphi,\psi]} \ind(f,g;p_1\!\Coin(\alpha\tilde{f},\tilde{g}))\,[\alpha] \;\in \ZZ\R[\varphi,\psi].
\end{align*}
A trace formula similar to \eqref{eq:L(f)} exists for the Lefschetz coincidence number, but in contrast with the fixed point case, not for the Reidemeister trace.

These invariants can also be generalised to the setting of $n$-valued maps, where an $n$-valued map $f:X\multimap X$ is a set-valued function such that $f(x)\subseteq X$ has cardinality exactly $n$ for each $x\in X$, and which is upper- and lower-semicontinuous in the sense of \cite{browngoncalves}. Equivalently, an $n$-valued map of $X$ is a single-valued map $f:X\to D_n(X)$, where \[
D_n(X)=\{\{x_1,\ldots,x_n\}\subseteq X \mid x_i\neq x_j \text{ if }i\neq j \}
\]
is the $n$-th unordered configuration space of $X$, topologised as the quotient of the ordered configuration space \[
F_n(X)=\{(x_1,\ldots,x_n) \mid x_i\neq x_j \text{ if }i\neq j \}\subseteq X^n
\]
by the action of the symmetric group, which we write as $\S_n$ (see \cite{browngoncalves}).

We say $x\in X$ is a fixed point of an $n$-valued map $f:X\to D_n(X)$ if $x\in f(x)$, and we write the set of fixed points of $f$ as $\Fix(f)$. Nielsen classes and the fixed point index were generalised to $n$-valued maps on compact polyhedra by Schirmer in \cite{schirmer2}. Reidemeister classes were generalised to this setting in \cite{charlotte1}. Let us briefly recall the main results regarding the latter.

Let $X$ be a compact polyhedron with universal cover $p:\tilde{X}\to X$ and covering group $\pi$. The $n$-th orbit configuration space with respect to $p$ is defined as \[
F_n(\tilde{X},\pi)=\{(\tilde{x}_1,\ldots,\tilde{x}_n) \mid p(\tilde{x}_i)\neq p(\tilde{x}_j) \text{ if }i\neq j \}\subseteq \tilde{X}^n,
\]
and the map \[
p^n:F_n(\tilde{X},\pi)\to D_n(X):(\tilde{x}_1,\ldots,\tilde{x}_n)\to \{p(\tilde{x}_1),\ldots,p(\tilde{x}_n) \}
\]
is a covering map. The group of covering transformations is isomorphic to $\pi^n\rtimes \S_n$, where the group operation is given by \[
(\alpha_1,\ldots,\alpha_n;\sigma)(\beta_1,\ldots,\beta_n;\tau)=(\alpha_1\beta_{\sigma^{-1}(1)},\ldots,\alpha_n\beta_{\sigma^{-1}(n)};\sigma\tau)
\] 
and where $(\alpha_1,\ldots,\alpha_n;\sigma)\in\pi^n\rtimes \S_n$ acts on $(\tilde{x}_1,\ldots,\tilde{x}_n)\in F_n(\tilde{X},\pi)$ by \[
(\alpha_1,\ldots,\alpha_n;\sigma)(\tilde{x}_1,\ldots,\tilde{x}_n)=(\alpha_1\tilde{x}_{\sigma^{-1}(1)},\ldots,\alpha_n\tilde{x}_{\sigma^{-1}(n)}).
\]

Given an $n$-valued map $f:X\to D_n(X)$, we can choose a lift $\tilde{f}:\tilde{X}\to F_n(\tilde{X},\pi)$ so that \[
\begin{tikzcd}
\tilde{X} \ar[r,"\tilde{f}"] \ar[d,"p"'] & F_n(\tilde{X},\pi) \ar[d,"p^n"] \\
X \ar[r,"f"] & D_n(X)
\end{tikzcd}
\]
commutes. As $F_n(\tilde{X},\pi)$ is a subspace of $\tilde{X}^n$, we can write $\tilde{f}=(\tilde{f}_1,\ldots,\tilde{f}_n)$ for single-valued maps $\tilde{f}_i:\tilde{X}\to\tilde{X}$.

This lift $\tilde{f}$ induces a morphism of the covering groups $\tilde{f}_\#:\pi\to \pi^n\rtimes \S_n$ satisfying \[
\tilde{f}\alpha=\tilde{f}_\#(\alpha)\tilde{f},\quad \forall \alpha\in \pi.
\]
If we write $\tilde{f}_\#=(\varphi_1,\ldots,\varphi_n;\sigma)$ with $\varphi_i:\pi\to \pi$ and $\sigma:\pi\to \S_n$, this is equivalent to \[
(\tilde{f}_1\alpha,\ldots,\tilde{f}_n\alpha)=(\varphi_1(\alpha)\tilde{f}_{\sigma_\alpha^{-1}(1)},\ldots,\varphi_n(\alpha)\tilde{f}_{\sigma_\alpha^{-1}(n)})
\]
where we write $\sigma_\alpha=\sigma(\alpha)\in \S_n$. The map $\sigma:\pi\to \S_n$ is a group morphism; the maps $\varphi_i:\pi\to \pi$ in general are not, but their restrictions to the subgroups $S_i=\{\alpha\in \pi\mid \sigma_\alpha(i)=i \}$ are.

\begin{rmk}\label{rmk:lift-choice}
As in the single-valued case, the morphism $\tilde{f}_\#$ depends on the choice of lift $\tilde{f}$: for any other lift $\tilde{f}'$, there is a covering translation $(\gamma_1,\ldots,\gamma_n;\eta)\in \pi^n\rtimes \S_n$ so that \[
\tilde{f}'=(\gamma_1,\ldots,\gamma_n;\eta)\tilde{f}=(\gamma_1\tilde{f}_{\eta^{-1}(1)},\ldots,\gamma_n\tilde{f}_{\eta^{-1}(n)}).
\]
The morphism induced by this lift is $\tilde{f}'_\#=(\varphi'_1,\ldots,\varphi'_n;\sigma')$, where $\sigma'_{\alpha}=\eta\sigma_{\alpha}\eta^{-1}$ and $\varphi'_i(\alpha)=\gamma_i\varphi_{\eta^{-1}(i)}(\alpha)\gamma^{-1}_{\eta\sigma_\alpha^{-1}\eta^{-1}(i)}$ for all $\alpha\in \pi$. In particular, the morphism obtained by restricting $\varphi'_i$ to $S'_i = \{\alpha\in \pi \mid \eta\sigma_\alpha\eta^{-1}(i)=i \}=S_{\eta^{-1}(i)}$ is $\tau_{\gamma_i}\varphi_{\eta^{-1}(i)}$.
\end{rmk}

The set of Reidemeister classes $\R[\tilde{f}_\#]$ is defined as the quotient of $\pi\times \{1,\ldots,n\}$ by the relation \[
(\alpha,i) \sim_{\tilde{f}_\#} (\beta,j) \iff \exists \gamma\in \pi: \sigma_\gamma(i)=j \text{ and } \alpha=\gamma\beta\varphi_i(\gamma^{-1}).
\]
The Reidemeister class of $(\alpha,i)$ will be denoted by $[(\alpha,i)]_{\tilde{f}_\#}$, or by $[(\alpha,i)]$ when $\tilde{f}_\#$ is understood.

\begin{rmk}\label{rmk:R(f)-R(phi,i)}
It follows from the reasoning made in \cite[Section 4]{charlotte1} that the Reidemeister classes $\R[\tilde{f}_\#]$ can be partitioned as follows: if $i_1,\ldots,i_r$ are representatives for the orbits of $\{1,\ldots,n\}$ under the action of $\pi$ given by $\gamma\cdot i=\sigma_\gamma(i)$, then \[
\R[\tilde{f}_\#]=\bigsqcup_{\ell=1}^r \bigsqcup_{[\alpha]\in \R[\varphi_{i_\ell},\iota_{i_\ell}]} [(\alpha,i_\ell)]
\]
where $\varphi_{i_\ell}$ is considered as a morphism $S_{i_\ell}\to \pi$ and $\iota_{i_\ell}:S_{i_\ell}\to \pi$ is the inclusion. Recall that $\R[\varphi_{i_\ell},\iota_{i_\ell}]$ denotes the set of Reidemeister coincidence classes of the morphisms $\varphi_{i_\ell}$ and $\iota_{i_\ell}$, as defined on page \pageref{R-coin}.

In particular, suppose we have to compute a sum of the form \[
\sum_{[(\alpha,i)]\in \R[\tilde{f}_\#]} X_{[(\alpha,i)]}
\]
where the terms $X_{[(\alpha,i)]}$ belong to an abelian group and only depend on the Reidemeister classes $[(\alpha,i)]$ of $(\alpha,i)$. Then this sum can be decomposed as \[
\sum_{[(\alpha,i)]\in \R[\tilde{f}_\#]} X_{[(\alpha,i)]} = \sum_{\ell=1}^r \sum_{[\alpha]\in \R[\varphi_{i_\ell},\iota_{i_\ell}]} X_{[(\alpha,i_\ell)]},
\] 
or also
\[
\sum_{[(\alpha,i)]\in \R[\tilde{f}_\#]} X_{[(\alpha,i)]} = \sum_{i=1}^n \frac{1}{[\pi:S_i]} \sum_{[\alpha]\in \R[\varphi_i,\iota_i]} X_{[(\alpha,i)]},
\]
since the partition is independent of the chosen representatives $i_\ell$, and the orbit of $i$ has $[\pi:S_i]$ elements, by the orbit-stabilizer theorem.
\end{rmk}

The Reidemeister classes stand in one-to-one correspondence with the Nielsen fixed point classes $p\Fix(\alpha\tilde{f}_i)$. Let $\ind(f;p\Fix(\alpha\tilde{f}_i))$ denote the index of a fixed point class, as defined by Schirmer in \cite{schirmer2}. Then the Lefschetz and Nielsen numbers of $f$ can be written as \begin{align*}
L(f)&=\sum_{[(\alpha,i)]\in \R[\tilde{f}_\#]} \ind(f;p\Fix(\alpha\tilde{f}_i)) \\
N(f)&=\#\{[(\alpha,i)]\in \R[\tilde{f}_\#] \mid \ind(f;p\Fix(\alpha\tilde{f}_i))\neq 0 \}.
\end{align*}
The Reidemeister trace of $f$ is defined as\footnote{In \cite{staecker2023}, another convention is used for labeling the Reidemeister classes, with the relation $\alpha=\varphi_j(\gamma)\beta\gamma^{-1}$ instead of $\alpha=\gamma\beta\varphi_i(\gamma^{-1})$, so that the Reidemeister class $[(\alpha,i)]$ represents the fixed point class $p\Fix(\alpha^{-1}\tilde{f}_i)$ instead of $p\Fix(\alpha\tilde{f}_i)$. This is merely a convention taken for technical reasons in \cite{staecker2023}; we opt for the more commonly used one, as e.g.\ in \cite{charlotte1}.} \[
\RT(f,\tilde{f})=\sum_{[(\alpha,i)]\in \R[\tilde{f}_\#]} \ind(f;p\Fix(\alpha\tilde{f}_i))\,[(\alpha,i)] \;\in \ZZ\R[\tilde{f}_\#].
\]

In \cite{staecker2023}, another averaging formula is shown:

\begin{thm}[{\cite[Corollary 27]{staecker2023}}]\label{thm:staecker}
Let $\Gamma$ be a finite index normal subgroup of $\pi$ so that $\tilde{f}_\#(\Gamma)\subseteq \Gamma^n\rtimes \S_n$, and let $\bar{X}=\Gamma\orb \tilde{X}$ be the corresponding finite cover of $X$. If $\bar{f}:\bar{X}\to D_n(\bar{X})$ is a lift of $f$, then \[
\RT(f,\tilde{f})=\frac{1}{[\pi:\Gamma]} \sum_{\bar{\alpha}\in \pi/\Gamma} \hat{r}^{\alpha}(\RT(\bar{\alpha}^n\bar{f},\alpha^n\tilde{f})),
\]
with $\hat{r}^\alpha:\ZZ\R[\tau_{\alpha^n}\bar{f}_\#]\to \ZZ\R[\tilde{f}_\#]$ given by $\hat{r}^\alpha([(\beta,i)]_{\tau_{\alpha^n}\bar{f}_\#})=[(\beta\alpha,i)]_{\bar{f}_\#}$ (here $\bar{\alpha}^n \bar{f}$ denotes the $n$-valued map given by $\bar{\alpha}^n \bar{f}(\bar{x})=\{\bar{\alpha}\bar{x}_1,\ldots,\bar{\alpha}\bar{x}_n\}$ if $\bar{f}(\bar{x})=\{\bar{x}_1,\ldots,\bar{x}_n\}$, and $\alpha^n\tilde{f}$ the one given by $\alpha^n \tilde{f}(\tilde{x})=(\alpha\tilde{x}_1,\ldots,\alpha\tilde{x}_n)$ if $\tilde{f}(\tilde{x})=(\tilde{x}_1,\ldots,\tilde{x}_n)$).
\end{thm}

However, the setting in which this theorem is applicable is much more restricted than in the single-valued case. For example, any single-valued map on an infra-nilmanifold (e.g.\ a compact flat manifold) admits a lift to a finite cover which is a nilmanifold (resp.\ a torus). In the $n$-valued case, this is no longer true, as demonstrated by the following example from \cite{affien}. For the convenience of the reader, we recall this example in full detail, because we will return to it later.

\begin{ex}\label{ex:klein}
Let $\pi\subseteq \RR^2\rtimes \mathrm{GL}_2(\RR)$ be the group generated by \[
a=\left(\begin{bmatrix}1 \\ 0 \end{bmatrix},\begin{bmatrix} 1 & 0 \\ 0 & 1 \end{bmatrix}\right), \quad b=\left(\begin{bmatrix}0 \\ \frac{1}{2} \end{bmatrix},\begin{bmatrix} -1 & 0 \\ 0 & 1 \end{bmatrix}\right).
\]
This group acts naturally on $\RR^2$ and the quotient $\pi\orb \RR^2$ is a compact flat manifold isomorphic to the Klein bottle. Consider the map given by \[
\tilde{f}: \RR^2\to F_2(\RR^2,\pi): \begin{bmatrix}t_1 \\ t_2 \end{bmatrix}\mapsto \left(\begin{bmatrix}0 \\ \frac{t_2}{2} \end{bmatrix},\begin{bmatrix}0 \\ \frac{t_2}{2}-\frac{1}{4} \end{bmatrix}\right).
\]
Note that for all $t_1,t_2\in \RR$
\begin{align*}
\tilde{f}\left(a\begin{bmatrix}t_1 \\ t_2 \end{bmatrix}\right)&=\tilde{f}\left(\begin{bmatrix}t_1+1 \\ t_2 \end{bmatrix}\right)
=\left(\begin{bmatrix}0 \\ \frac{t_2}{2} \end{bmatrix},\begin{bmatrix}0 \\ \frac{t_2}{2}-\frac{1}{4} \end{bmatrix}\right)=(1,1;1)\,\tilde{f}\left(\begin{bmatrix}t_1 \\ t_2 \end{bmatrix}\right) \\
\tilde{f}\left(b\begin{bmatrix}t_1 \\ t_2 \end{bmatrix}\right)&=\tilde{f}\left(\begin{bmatrix}-t_1 \\ t_2+\frac{1}{2} \end{bmatrix}\right)
=\left(\begin{bmatrix}0 \\ \frac{t_2}{2}+\frac{1}{4} \end{bmatrix},\begin{bmatrix}0 \\ \frac{t_2}{2} \end{bmatrix}\right)=(b,1;\sigma)\,\tilde{f}\left(\begin{bmatrix}t_1 \\ t_2 \end{bmatrix}\right)
\end{align*}
with $\sigma\in \S_2$ the permutation that interchanges $1$ and $2$. Thus, $\tilde{f}$ induces a well-defined $2$-valued map $f:\pi\orb \RR^2\to D_2(\pi\orb \RR^2)$ on the Klein bottle with morphism $\tilde{f}_\#:\pi\to \pi^2\rtimes \S_2$ given by $\tilde{f}_\#(a)=(1,1;1)$ and $\tilde{f}_\#(b)=(b,1;\sigma)$.\footnote{Note that in the notation $\tilde{f}_\#(a)=(1,1;1)$, the first two entries $1$ denote the identity element in $\pi$, i.e.\ $\left(\begin{bmatrix}0 \\ 0 \end{bmatrix},\begin{bmatrix} 1 & 0 \\ 0 & 1 \end{bmatrix}\right)$, and the $1$ after the semicolon denotes the identity permutation $1\in \S_2$.}

There is no finite index normal subgroup $\Gamma\subseteq \pi$ contained in $\RR^2$ so that $\tilde{f}_\#(\Gamma)\subseteq \Gamma^2\rtimes \S_2$ (in other words, there is no torus $\Gamma\orb \RR^2$ that finitely covers $\pi\orb \RR^2$ so that $f$ lifts to $\Gamma\orb \RR^2$). Indeed, note that $\tilde{f}_\#(b^2)=(b,b;1)$. If $\Gamma$ is a finite index subgroup of $\pi$ contained in $\pi\cap \RR^2=\langle a,b^2 \rangle$, there is a smallest $\ell\in\ZZ_{>0}$ so that $(b^2)^\ell\in \Gamma$. But then $\tilde{f}_\#((b^2)^\ell)=(b^\ell,b^\ell;1)$ is not an element of $\Gamma^2\rtimes \S_2$, by minimality of $\ell$.
\end{ex}

In this paper, we show another type of averaging formula for the Reidemeister trace, in terms of coincidence traces of single-valued maps. We lose a bit of generality since, by passing onto coincidences, we need our cover to be a closed orientable manifold, and hence our space itself to be a closed manifold. On the other hand, in classical settings like for infra-nilmanifolds, our type of lift will always exist (unlike the one from Theorem \ref{thm:staecker}), making our formula more widely applicable.

\section{The split lift of an $n$-valued map}\label{sec:split-lift}

Let $X$ be a compact polyhedron with universal cover $p:\tilde{X}\to X$ and group of covering transformations $\pi$. Let $\Gamma$ be a finite index normal subgroup of $\pi$, and $\bar{X}=\Gamma\orb \tilde{X}$ the corresponding finite cover of $X$, with covering map $\bar{p}:\bar{X}\to X$. Let $f:X\to D_n(X)$ be an $n$-valued map with lift $\tilde{f}=(\tilde{f}_1,\ldots,\tilde{f}_n):\tilde{X}\to F_n(\tilde{X},\pi)$ and induced morphism $\tilde{f}_\#=(\varphi_1,\ldots,\varphi_n;\sigma):\pi\to \pi^n\rtimes \S_n$.

\begin{df}
We will say a finite index normal subgroup $S$ of $\pi$ is \emph{$f$-$\Gamma$-invariant} if $S$ is contained in $\Gamma$ and in $\ker\sigma$ (hence in all groups $S_i=\{\alpha\in \pi \mid \sigma_\alpha(i)=i \}$), and $\varphi_i(S)\subseteq \Gamma$ for all $i$.
\end{df}

Note that this definition is independent of the chosen lift $\tilde{f}$: as we saw in Remark \ref{rmk:lift-choice}, another lift induces morphisms $\sigma'=\eta\sigma\eta^{-1}$, with the same kernel as $\sigma$, and $\varphi'_i=\tau_{\gamma_i}\varphi_{\eta^{-1}(i)}$, where $\tau_{\gamma_i}\varphi_{\eta^{-1}(i)}(S)\subseteq \Gamma$ if and only if $\varphi_{\eta^{-1}(i)}(S)\subseteq \tau_{\gamma_i}^{-1}(\Gamma)=\Gamma$.

\begin{rmk}\label{rmk:virt-polyc}
If $\pi$ is virtually polycyclic, such group $S$ always exists: for instance, \[
S=\langle \alpha^{[\pi:\Gamma]} \mid \alpha\in \ker\sigma \rangle.
\]
Clearly, this is a normal subgroup of $\pi$ contained in $\Gamma$ and in $\ker\sigma$, and $\varphi(S)\subseteq \Gamma$ for any morphism $\varphi:\ker\sigma\to \pi$. Since $\pi/S$ is periodic, and periodic virtually polycyclic groups are automatically finite, the group $S$ is of finite index in $\pi$.
\end{rmk}

Suppose $S$ is an $f$-$\Gamma$-invariant subgroup of $\pi$. Then we can consider $\hat{X}=S\orb \tilde{X}$ as a finite cover of $X$; let $\hat{p}:\hat{X}\to X$ be the covering map. The projection $p':\tilde{X}\to \hat{X}$ is the universal cover of $\hat{X}$, and satisfies $\hat{p}\circ p'=p$.  We also have a finite covering map $q:\hat{X}\to \bar{X}$ so that $\bar{p}\circ q=\hat{p}$, and $q'=q\circ p':\tilde{X}\to \bar{X}$ is the universal cover of $\bar{X}$. Together: \[
\begin{tikzcd}
\tilde{X} \ar[r,"p'"] \ar[drr,"p"',yshift=-1pt] \ar[rr,bend left=20,yshift=4pt,"q'"] & \hat{X} \ar[r,"q"] \arrow[xshift=1pt,yshift=2pt]{dr}[description]{\hat{p}} & \bar{X} \ar[d,"\bar{p}"] \\
& & X.
\end{tikzcd}
\]

As a finite covering space for $D_n(X)$, consider \[
F_n(\bar{X},\pi/\Gamma)=\{(\bar{x}_1,\ldots,\bar{x}_n)\in \bar{X}^n \mid \bar{p}(\bar{x}_i)\neq \bar{p}(\bar{x}_j) \text{ if }i\neq j \}
\] 
with covering map \[
\bar{p}^n:F_n(\bar{X},\pi/\Gamma)\to D_n(X):(\bar{x}_1,\ldots,\bar{x}_n)\mapsto \{\bar{p}(\bar{x}_1),\ldots,\bar{p}(\bar{x}_n) \}.
\]
The map $(q')^n:F_n(\tilde{X},\pi)\to F_n(\bar{X},\pi/\Gamma)$ obtained by applying the universal cover $q'$ of $\bar{X}$ in each coordinate is also a well-defined covering map.

We obtain a commutative diagram \[
\begin{tikzcd}
\tilde{X} \ar[d,"p'"'] \ar[r,"\tilde{f}"] & F_n(\tilde{X},\pi) \ar[d,"(q')^n"] \\
\hat{X} \ar[d,"\hat{p}"'] & F_n(\bar{X},\pi/\Gamma) \ar[d,"\bar{p}^n"] \\
X \ar[r,"f"] & D_n(X).
\end{tikzcd}
\]

\begin{lem}
There is a (unique) map $\bar{f}:\hat{X}\to F_n(\bar{X},\pi/\Gamma)$ making the above diagram commutative.
\end{lem}

\begin{proof}
If such a map $\bar{f}$ exists, it must be given by $\bar{f}(p'(\tilde{x}))=(q')^n(\tilde{f}(\tilde{x}))$; so it suffices to show that this map is well-defined. (Indeed, if $\bar{f}$ is given by this formula, also the lower half of the diagram commutes, since the big diagram commutes and $p'$ is surjective.) 

Suppose $p'(\tilde{x})=p'(\tilde{x}')$; then there is a $\gamma\in S$ such that $\tilde{x}'=\gamma\tilde{x}$, according to the action of $S$ on $\tilde{X}$. We have to show that $(q')^n(\tilde{f}(\gamma\tilde{x}))=(q')^n(\tilde{f}(\tilde{x}))$. 

By definition of the maps $\varphi_i$, and since $S\subseteq \ker\sigma$, we have
\begin{align*}
(q')^n(\tilde{f}(\gamma\tilde{x}))&=(q')^n(\tilde{f}_1(\gamma\tilde{x}),\ldots,\tilde{f}_n(\gamma\tilde{x})) \\ &=(q')^n(\varphi_1(\gamma)\tilde{f}_1(\tilde{x}),\ldots,\varphi_n(\gamma)\tilde{f}_n(\tilde{x})) \\
&=(q'(\varphi_1(\gamma)\tilde{f}_1(\tilde{x})),\ldots,q'(\varphi_n(\gamma)\tilde{f}_n(\tilde{x}))).
\end{align*}
Since $\varphi_i(S)\subseteq \Gamma$ for each $i$, and $\Gamma$ is the covering group of $q'$, the last line is equal to \[
(q'(\tilde{f}_1(\tilde{x})),\ldots,q'(\tilde{f}_n(\tilde{x})))=(q')^n(\tilde{f}(\tilde{x})). \qedhere
\]
\end{proof}

Since $F_n(\bar{X},\pi/\Gamma)$ is a subspace of $\bar{X}^n$, there are single-valued maps $\bar{f}_i:\hat{X}\to \bar{X}$ so that $\bar{f}=(\bar{f}_1,\ldots,\bar{f}_n)$. Note that, by construction, $\tilde{f}_i:\tilde{X}\to \tilde{X}$ is a lift of $\bar{f}_i$ for each $i$.

\begin{df}
We call the map $\bar{f}=(\bar{f}_1,\ldots,\bar{f}_n)$ the \emph{split lift of $f$ with respect to $S$ and $\Gamma$, corresponding to $\tilde{f}$}.
\end{df}

\begin{rmk}
A point $x=\hat{p}(\hat{x})\in X$ is a fixed point of $f$ if and only if \begin{align*}
\hat{p}(\hat{x})\in f(\hat{p}(\hat{x})) &\iff \hat{p}(\hat{x})\in \bar{p}^n(\bar{f}(\hat{x}))=\{\bar{p}(\bar{f}_1(\hat{x})),\ldots,\bar{p}(\bar{f}_n(\hat{x})) \} \\&\iff \exists i\in \{1,\ldots,n\}: \hat{p}(\hat{x})=\bar{p}(\bar{f}_i(\hat{x})) \\&\iff \exists i\in \{1,\ldots,n\}: \bar{p}(q(\hat{x}))=\bar{p}(\bar{f}_i(\hat{x})) \\&\iff \exists i\in \{1,\ldots,n\}, \exists\bar{\alpha}\in \pi/\Gamma: q(\hat{x})=\bar{\alpha}\bar{f}_i(\hat{x}).
\end{align*}
That is, \[
\Fix(f)=\bigcup_{i=1}^n\bigcup_{\bar{\alpha}\in \pi/\Gamma}\hat{p}\Coin(\bar{\alpha}\bar{f}_i,q).
\]
We can thus relate fixed points of the $n$-valued map $f$ with coincidences of the single-valued maps $\bar{\alpha}\bar{f}_i$ and the projection $q$.
\end{rmk}

As reference lifts for $\bar{\alpha}\bar{f}_i:\hat{X}\to \bar{X}$ and $q:\hat{X}\to \bar{X}$ with respect to the universal coverings $p':\tilde{X}\to\hat{X}$ and $q':\tilde{X}\to \bar{X}$, we may respectively choose the maps $\alpha\tilde{f}_i:\tilde{X}\to \tilde{X}$ (for some representative $\alpha\in \pi$ of $\bar{\alpha}$) and the identity map $\id:\tilde{X}\to \tilde{X}$:
\[
\begin{tikzcd}
\tilde{X} \ar[r,"\alpha\tilde{f}_i","\id"'] \ar[d,"p'"'] & \tilde{X} \ar[d,"q'"] \\
\hat{X} \ar[r,"\bar{\alpha}\bar{f}_i","q"'] & \bar{X}.
\end{tikzcd}
\]
The morphism of covering groups induced by the lift $\alpha\tilde{f}_i$ of $\bar{\alpha}\bar{f}_i$ is the restriction $\tau_\alpha\varphi_i':S\to \Gamma$: for all $\gamma\in S$, \[
\alpha\tilde{f}_i\gamma=\alpha\varphi_i(\gamma)\tilde{f}_i=\tau_\alpha(\varphi_i(\gamma))\alpha\tilde{f}_i.
\]
The morphism induced by the lift $\id$ of $q$ is the inclusion morphism $\iota:S\to \Gamma$.

\section{Averaging formula for the Reidemeister trace}

We stay in the setting of the previous section, but from now on, we assume $X$ is a closed manifold and the finite cover $\bar{X}$ is a closed orientable manifold (e.g.\ the double orientable cover of $X$). Note that $\hat{X}$, being a finite cover of $\bar{X}$, is also a closed orientable manifold of the same dimension.

Recall that the Reidemeister trace of our $n$-valued map $f$ is given by \begin{equation}\label{eq:RT-nv}
\RT(f,\tilde{f})=\sum_{[(\alpha,i)]\in \R[\tilde{f}_\#]} \ind(f;p\Fix(\alpha\tilde{f}_i))\,[(\alpha,i)] \;\in \ZZ\R[\tilde{f}_\#].
\end{equation}
Since $\bar{\alpha}\bar{f}_i,q:\hat{X}\to \bar{X}$ are maps between equal-dimensional closed orientable manifolds, we can consider the coincidence Reidemeister trace: \[
\RT(\bar{\alpha}\bar{f}_i,\alpha\tilde{f}_i;q,\id)=\sum_{[\beta]\in \R[\tau_\alpha\varphi'_i,\iota]} \ind(\bar{\alpha}\bar{f}_i,q;p'\!\Coin(\beta\alpha\tilde{f}_i,\id))\,[\beta] \;\in \ZZ\R[\tau_\alpha\varphi'_i,\iota].
\]

For all $\alpha$ and $i$, we have a well-defined map \[
\hat{r}_i^\alpha: \R[\tau_\alpha\varphi'_i,\iota]\to \R[\tilde{f}_\#]: [\beta]_{\tau_\alpha\varphi'_i,\iota}\mapsto [(\beta\alpha,i)]_{\tilde{f}_\#}.
\]
Indeed, if $\beta\sim_{\tau_\alpha\varphi'_i,\iota}\beta'$ in $\pi$, then \[
\exists \gamma\in S: \beta=\gamma\beta'\alpha\varphi_i(\gamma^{-1})\alpha^{-1}.
\]
Since $S\subseteq S_i=\{\gamma\in \pi\mid \sigma_\gamma(i)=i \}$, this implies \[
\exists \gamma\in \pi: \sigma_\gamma(i)=i \text{ and } \beta\alpha=\gamma\beta'\alpha\varphi_i(\gamma^{-1}),
\]
which means $(\beta\alpha,i)\sim_{\tilde{f}_\#}(\beta'\alpha,i)$.

We can consider the extension of $\hat{r}_i^\alpha$ to $\ZZ \R[\tau_\alpha\varphi'_i,\iota]\to \ZZ \R[\tilde{f}_\#]$, which we will still write as $\hat{r}_i^\alpha$. With these notations, we will show:

\begin{thm}[Averaging formula]\label{thm:av-RT}
For an $n$-valued map $f:X\to D_n(X)$ with lift $\tilde{f}=(\tilde{f}_1,\ldots,\tilde{f}_n):\tilde{X}\to F_n(\tilde{X},\pi)$ and corresponding split lift $(\bar{f}_1,\ldots,\bar{f}_n)$ with respect to $S$ and $\Gamma$,
\[
\RT(f,\tilde{f})=\frac{1}{[\pi:S]}\sum_{i=1}^{n}\sum_{\bar{\alpha}\in \pi/\Gamma} \hat{r}_i^\alpha(\RT(\bar{\alpha}\bar{f}_i,\alpha\tilde{f}_i;q,\id)).
\]
\end{thm}

First of all, as observed in Remark \ref{rmk:R(f)-R(phi,i)}, we can decompose a sum over the Reidemeister classes of $\tilde{f}_\#$ as \begin{equation}\label{eq:decomp1}
\sum_{[(\alpha,i)]\in \R[\tilde{f}_\#]} X_{[(\alpha,i)]} = \sum_{i=1}^n \frac{1}{[\pi:S_i]} \sum_{[\alpha]\in \R[\varphi_i,\iota_i]} X_{[(\alpha,i)]}
\end{equation}
where $\varphi_i$ should be considered as a morphism $S_i\to \pi$, and $\iota_i:S_i\to \pi$ is the inclusion. The second sum can be decomposed further, using the technique we developed in \cite[Section 3]{affien}. For the convenience of the reader, we show how the reasoning of \cite{affien} can be adapted to this setting.

\begin{lem}\label{lem:decomp}
Consider a sum of the type \[
\sum_{[(\alpha,i)]\in \R[\tilde{f}_\#]} X_{[(\alpha,i)]}
\]
where $X_{[(\alpha,i)]}$ are elements of an abelian group, indexed by the Reidemeister classes of $\tilde{f}_\#$. This sum can be decomposed as \[
\sum_{[(\alpha,i)]\in \R[\tilde{f}_\#]} X_{[(\alpha,i)]} = \frac{1}{[\pi:S]} \sum_{i=1}^n \sum_{\bar{\alpha}\in \pi/\Gamma} \sum_{[\beta]\in \R[\tau_\alpha\varphi'_i,\iota]} |u_i(\coin(\tau_{\beta\alpha}\varphi_i,\iota_i))|\,X_{[(\beta\alpha,i)]}
\]
where $u_i:S_i\to S_i/S$ is the projection, $\coin(\tau_{\beta\alpha}\varphi_i,\iota_i)\subseteq S_i$ is the coincidence subgroup of the morphisms $\tau_{\beta\alpha}\varphi_i,\iota_i:S_i\to \pi$, and the bars denote cardinality.
\end{lem}

In the proof of this lemma, we will use the following result from \cite{HLP2012}:

\begin{lem}[{\cite[{Lemma 3.1.4}]{HLP2012}}]\label{lem:HLP}
Consider a commutative diagram of groups and group morphisms with exact rows \[
\begin{tikzcd}
1 \ar[r] & \Gamma_1 \ar[r,"i_1"] \ar[d,"\varphi'"',"\psi'"] & \Pi_1 \ar[r,"u_1"] \ar[d,"\varphi"',"\psi"] & \Pi_1/\Gamma_1 \ar[r] \ar[d,"\bar{\varphi}"',"\bar{\psi}"] & 1 \\
1 \ar[r] & \Gamma_2 \ar[r,"i_2"] & \Pi_2 \ar[r,"u_2"] & \Pi_2/\Gamma_2 \ar[r] & 1
\end{tikzcd}
\]
such that the quotient groups $\Pi_i/\Gamma_i$ are finite. For all $\alpha\in \Pi_2$ and $\beta\in \Gamma_2$, the number of Reidemeister coincidence classes $[\beta']\in \R[\tau_\alpha\varphi',\psi']$ such that $\beta\sim_{\tau_\alpha\varphi,\psi}\beta'$ is equal to \[
[\coin(\tau_{\bar{\alpha}}\bar{\varphi},\bar{\psi}):u_1(\coin(\tau_{\beta\alpha}\varphi,\psi))].
\]
\end{lem}

\begin{proof}[Proof of Lemma \ref{lem:decomp}]
To decompose the second sum \[
\sum_{[\alpha]\in \R[\varphi_i,\iota_i]} X_{[(\alpha,i)]}
\] 
in \eqref{eq:decomp1}, note that we can write \[
\R[\varphi_i,\iota_i]=\bigcup_{\bar{\alpha}\in \pi/\Gamma}\,\bigcup_{\beta\in \Gamma}\,[\beta\alpha].
\]
If $\beta\alpha\sim_{\varphi_i,\iota_i}\beta'\alpha'$, then \[
\exists \gamma\in S_i: \beta'\alpha' = \gamma \beta\alpha \varphi_i(\gamma^{-1}).
\]
In particular, \[
\exists \bar{\gamma}\in S_i/S: \bar{\alpha}' = \bar{\iota}(\bar{\gamma}) \bar{\alpha} \bar{\varphi}_i(\bar{\gamma}^{-1})
\]
where $\bar{\iota}_i,\bar{\varphi}_i:S_i/S\to \pi/\Gamma$ are the morphisms induced by $\iota_i,\varphi_i:S_i\to \pi$;
that is, $\bar{\alpha}\sim_{\bar{\varphi}_i,\bar{\iota}_i}\bar{\alpha}'$. Conversely, if $\bar{\alpha}\sim_{\bar{\varphi}_i,\bar{\iota}_i}\bar{\alpha}'$, then there is a $\beta\in \Gamma$ so that $\beta\alpha\sim_{\varphi_i,\iota_i}\alpha'$.

For fixed $\alpha\in \pi$, we have $\beta\alpha\sim_{\varphi_i,\iota_i}\beta'\alpha$ if and only if \[
\exists \gamma\in S_i:\beta=\gamma \beta'\alpha\varphi_i(\gamma^{-1})\alpha^{-1}=\gamma\beta'\tau_\alpha\varphi_i(\gamma^{-1});
\] 
that is, $\beta\sim_{\tau_\alpha\varphi_i,\iota_i}\beta'$.
Thus, the above are disjoint unions when taken over these respective Reidemeister classes: \[
\R[\varphi_i,\iota_i]=\bigsqcup_{[\bar{\alpha}]\in \R[\bar{\varphi}_i,\bar{\iota}_i]}\,\bigsqcup_{[\beta]\in \R[\tau_\alpha\varphi_i,\iota_i],\beta\in \Gamma}\,[\beta\alpha].
\]

Rather than over the elements of $\R[\tau_\alpha\varphi_i,\iota_i]$ with representatives in $\Gamma$, we would like the second union to go over $[\beta]\in \R[\tau_\alpha\varphi'_i,\iota]$, where \[
\beta\sim_{\tau_\alpha\varphi'_i,\iota} \beta' \iff \exists \gamma\in S: \beta=\gamma \beta' \tau_\alpha\varphi_i(\gamma^{-1}).
\]
Clearly, $\beta\sim_{\tau_\alpha\varphi'_i,\iota} \beta'$ implies $\beta\sim_{\tau_\alpha\varphi_i,\iota_i} \beta'$. For the converse, we use Lemma \ref{lem:HLP} abovem which describes how many different classes in $\R[\tau_\alpha\varphi'_i,\iota]$ correspond to the same class in $\R[\tau_\alpha\varphi_i,\iota_i]$. Applying this to the diagram \[
\begin{tikzcd}
1 \ar[r] & S \ar[r] \ar[d,"\varphi'_i"',"\iota"] & S_i \ar[r,"u_i"] \ar[d,"\varphi_i"',"\iota_i"] & S_i/S \ar[r] \ar[d,"\bar{\varphi}_i"',"\bar{\iota}_i"] & 1 \\
1 \ar[r] & \Gamma \ar[r] & \pi \ar[r] & \pi/\Gamma \ar[r] & 1
\end{tikzcd}
\]
shows that for a given $\beta\in\Gamma$, the number of classes $[\beta']\in \R[\tau_\alpha\varphi'_i,\iota]$ such that $\beta\sim_{\tau_\alpha\varphi_i,\iota_i}\beta'$ is equal to \[
[\coin(\tau_{\bar{\alpha}}\bar{\varphi}_i,\bar{\iota}_i):u_i(\coin(\tau_{\beta\alpha}\varphi_i,\iota_i))]=\frac{|\coin(\tau_{\bar{\alpha}}\bar{\varphi}_i,\bar{\iota}_i)|}{|u_i(\coin(\tau_{\beta\alpha}\varphi_i,\iota_i))|}.
\]
Thus, we can write \[
\sum_{[\alpha]\in \R[\varphi_i,\iota_i]} X_{[(\alpha,i)]}=\sum_{[\bar{\alpha}]\in \R[\bar{\varphi}_i,\bar{\iota}_i]}\,\sum_{[\beta]\in \R[\tau_\alpha\varphi'_i,\iota]} \frac{|u_i(\coin(\tau_{\beta\alpha}\varphi_i,\iota_i))|}{|\coin(\tau_{\bar{\alpha}}\bar{\varphi}_i,\bar{\iota}_i)|} X_{[(\beta\alpha,i)]}.
\]

Lastly, consider the action of $S_i/S$ on $\pi/\Gamma$ given by \[
\bar{\gamma}\cdot \bar{\alpha}=\bar{\iota}_i(\bar{\gamma})\bar{\alpha}\bar{\varphi}_i(\bar{\gamma}^{-1}).
\]
The orbit of an element $\bar{\alpha}$ for this action is its Reidemeister class $[\bar{\alpha}]_{\bar{\varphi}_i,\bar{\iota}_i}$. The stabiliser of $\bar{\alpha}$ is the set \[
\{\bar{\gamma}\in S_i/S \mid \bar{\iota}_i(\bar{\gamma})\bar{\alpha}\bar{\varphi}_i(\bar{\gamma}^{-1})=\bar{\alpha} \}=\coin(\tau_{\bar{\alpha}}\bar{\varphi}_i,\bar{\iota}_i).
\]
Therefore we have \[
|[\bar{\alpha}]_{\bar{\varphi}_i,\bar{\iota}_i}|\cdot |\coin(\tau_{\bar{\alpha}}\bar{\varphi}_i,\bar{\iota}_i)|=|S_i/S|=[S_i:S]
\]
and we can rewrite the above sum as \begin{align*}
\sum_{[\alpha]\in \R[\varphi_i,\iota_i]} X_{[(\alpha,i)]} &= \sum_{\bar{\alpha}\in \pi/\Gamma}\frac{1}{|[\bar{\alpha}]_{\bar{\varphi}_i,\bar{\iota}_i}|}\sum_{[\beta]\in \R[\tau_\alpha\varphi'_i,\iota]} \frac{|u_i(\coin(\tau_{\beta\alpha}\varphi_i,\iota_i))|}{|\coin(\tau_{\bar{\alpha}}\bar{\varphi}_i,\bar{\iota}_i)|} X_{[(\beta\alpha,i)]} \\
&= \frac{1}{[S_i:S]}\sum_{\bar{\alpha}\in \pi/\Gamma}\sum_{[\beta]\in \R[\tau_\alpha\varphi'_i,\iota]} |u_i(\coin(\tau_{\beta\alpha}\varphi_i,\iota_i))| X_{[(\beta\alpha,i)]}.
\end{align*}
Together, we get \begin{align*}
\sum_{[(\alpha,i)]\in \R[\tilde{f}_\#]} X_{[(\alpha,i)]} &= \sum_{i=1}^n \frac{1}{[\pi:S_i]} \sum_{[\alpha]\in \R[\varphi_i,\iota_i]} X_{[(\alpha,i)]} \\
&= \sum_{i=1}^n \frac{1}{[\pi:S_i]} \frac{1}{[S_i:S]}\sum_{\bar{\alpha}\in \pi/\Gamma}\sum_{[\beta]\in \R[\tau_\alpha\varphi'_i,\iota]} |u_i(\coin(\tau_{\beta\alpha}\varphi,\psi))| X_{[(\beta\alpha,i)]} \\
&= \frac{1}{[\pi:S]} \sum_{i=1}^n \sum_{\bar{\alpha}\in \pi/\Gamma} \sum_{[\beta]\in \R[\tau_\alpha\varphi'_i,\iota]} |u_i(\coin(\tau_{\beta\alpha}\varphi_i,\iota_i))|\,X_{[(\beta\alpha,i)]}. \qedhere
\end{align*}
\end{proof}

By applying Lemma \ref{lem:decomp} to the terms \[
X_{[(\alpha,i)]}=\ind(f;p\Fix(\alpha\tilde{f}_i))\,[(\alpha,i)]_{\tilde{f}_\#}
\]
in the abelian group $\ZZ\R[\tilde{f}_\#]$, the formula \eqref{eq:RT-nv} for $\RT(f,\tilde{f})$ can be rewritten as \[
\frac{1}{[\pi:S]}\sum_{i=1}^{n}\sum_{\bar{\alpha}\in \pi/\Gamma}\sum_{[\beta]\in \R[\tau_\alpha\varphi'_i,\iota]} |u_i(\coin(\tau_{\beta\alpha}\varphi_i,\iota_i))|\ind(f;p\Fix(\beta\alpha\tilde{f}_i))\,[(\beta\alpha,i)]_{\tilde{f}_\#}
\] 
or, using the map $\hat{r}_i^\alpha$, as \[
\frac{1}{[\pi:S]}\sum_{i=1}^{n}\sum_{\bar{\alpha}\in \pi/\Gamma}\sum_{[\beta]\in \R[\tau_\alpha\varphi'_i,\iota]} |u_i(\coin(\tau_{\beta\alpha}\varphi_i,\iota_i))|\ind(f;p\Fix(\beta\alpha\tilde{f}_i))\;\hat{r}_i^\alpha([\beta]_{\tau_\alpha\varphi'_i,\iota}).
\]

To prove Theorem \ref{thm:av-RT}, it remains to show that \[
|u_i(\coin(\tau_{\beta\alpha}\varphi_i,\iota_i))|\,\ind(f;p\Fix(\beta\alpha\tilde{f}_i))=\ind(\bar{\alpha}\bar{f}_i,q;p'\!\Coin(\beta\alpha\tilde{f}_i,\id))
\]
for all $\alpha\in \pi$, $\beta\in \Gamma$ and $i\in \{1,\ldots,n\}$. For this, we will use the following auxiliary result.

\begin{lem}
For fixed $\alpha\in \pi$ and $i\in \{1,\ldots,n\}$, for all $x\in p\Fix(\alpha\tilde{f}_i)$, \[
|p'\Fix(\alpha\tilde{f}_i)\cap \hat{p}^{-1}(x)|=|u_i(\coin(\tau_\alpha\varphi_i,\iota_i))|.
\]
\end{lem}

\begin{proof}
Fix an element $\hat{x}\in p'\Fix(\alpha\tilde{f}_i)\cap \hat{p}^{-1}(x)$, and take $\tilde{x}\in \tilde{X}$ so that $\hat{x}=p'(\tilde{x})$ and $\alpha\tilde{f}_i(\tilde{x})=\tilde{x}$. We will express all elements of $p'\Fix(\alpha\tilde{f}_i)\cap \hat{p}^{-1}(x)$ in terms of $\hat{x}$.
Any other element in $\hat{p}^{-1}(x)$ is of the form $\bar{\gamma}\hat{x}$ for some $\bar{\gamma}\in \pi/S$. This element $\bar{\gamma}\hat{x}$ lies in $p'\Fix(\alpha\tilde{f}_i)$ if and only if it is of the form $p'(\tilde{x}')$ with $\alpha\tilde{f}_i(\tilde{x}')=\tilde{x}'$. Since $p(\tilde{x}')=x=p(\tilde{x})$, there is a $\gamma\in \pi$ such that $\tilde{x}'=\gamma\tilde{x}$. By covering space theory, it follows from $p'(\gamma\tilde{x})=\bar{\gamma}\hat{x}=\bar{\gamma}p'(\tilde{x})$ that $\bar{\gamma}$ is the image of $\gamma$ under the projection $\pi\to \pi/S$.

Since $\alpha\tilde{f}_i(\tilde{x}')=\tilde{x}'$, the element $\gamma$ satisfies \[
\alpha\tilde{f}_i(\gamma\tilde{x})=\gamma\tilde{x}=\gamma\alpha\tilde{f}_i(\tilde{x}) \iff \tilde{f}_i(\gamma\tilde{x})=\alpha^{-1}\gamma\alpha\tilde{f}_i(\tilde{x}).
\]
It follows that $\sigma_\gamma(i)=i$ and $\varphi_i(\gamma)=\alpha^{-1}\gamma\alpha$. In other words, $\gamma\in S_i$ and $\tau_\alpha\varphi_i(\gamma)=\gamma$; together, $\gamma\in \coin(\tau_\alpha\varphi_i,\iota_i)$.

We can conclude that \begin{align*}
p'\Fix(\alpha\tilde{f}_i)\cap \hat{p}^{-1}(x)&=\{\bar{\gamma}\hat{x}\mid \exists\gamma\in \coin(\tau_\alpha\varphi_i,\iota_i):\bar{\gamma}=u_i(\gamma) \} \\
&=\{\bar{\gamma}\hat{x}\mid \bar{\gamma}\in u_i(\coin(\tau_\alpha\varphi_i,\iota_i)) \}.
\end{align*}
In particular, $|p'\Fix(\alpha\tilde{f}_i)\cap \hat{p}^{-1}(x)|=|u_i(\coin(\tau_\alpha\varphi_i,\iota_i))|$.
\end{proof}

\begin{lem}\label{lem:ind}
For all $\alpha\in \pi$, $\beta\in \Gamma$ and $i\in \{1,\ldots,n\}$, \[
|u_i(\coin(\tau_{\beta\alpha}\varphi_i,\iota_i))|\,\ind(f;p\Fix(\beta\alpha\tilde{f}_i))=\ind(\bar{\alpha}\bar{f}_i,q;p'\!\Coin(\beta\alpha\tilde{f}_i,\id)).
\]
\end{lem}

\begin{proof}
After applying a homotopy if necessary, we may assume $f$ has only finitely many fixed points and \[
\ind(f;p\Fix(\beta\alpha\tilde{f}_i))=\sum_{x\in p\Fix(\beta\alpha\tilde{f}_i)} \ind(f;x)
\]
where $\ind(f;x)=\ind(\beta\alpha\tilde{f}_i;\tilde{x})$ for any choice of
$\tilde{x}\in\Fix(\beta\alpha\tilde{f}_i)\cap p^{-1}(x)$ (see \cite{schirmer2}). In the same spirit, \[
\ind(\bar{\alpha}\bar{f}_i,q;p'\!\Coin(\beta\alpha\tilde{f}_i,\id))=\sum_{\hat{x}\in p'\!\Coin(\beta\alpha\tilde{f}_i,\id)} \ind(\bar{\alpha}\bar{f}_i,q;\hat{x})
\]
where $\ind(\bar{\alpha}\bar{f}_i,q;\hat{x})=\ind(\beta\alpha\tilde{f}_i,\id;\tilde{x})$ for any $\tilde{x}\in\Coin(\beta\alpha\tilde{f}_i,\id)\cap (p')^{-1}(\hat{x})$. 

Note that $\Coin(\beta\alpha\tilde{f}_i,\id)=\Fix(\beta\alpha\tilde{f}_i)$ and, for all $\hat{x}\in p'\Fix(\beta\alpha\tilde{f}_i)$ and for any $\tilde{x}\in \Fix(\beta\alpha\tilde{f}_i)\cap (p')^{-1}(\hat{x})$, \begin{align*}
\ind(\bar{\alpha}\bar{f}_i,q;\hat{x})&=\ind(\beta\alpha\tilde{f}_i,\id;\tilde{x})\\ &=\ind(\beta\alpha\tilde{f}_i;\tilde{x})\\ &=\ind(f;p(\tilde{x}))\\ &=\ind(f;\hat{p}(\hat{x})).
\end{align*}
It follows that \begin{align*}
\ind(\bar{\alpha}\bar{f}_i,q;p'\!\Coin(\beta\alpha\tilde{f}_i,\id))&=\sum_{\hat{x}\in p'\Fix(\beta\alpha\tilde{f}_i)} \ind(f;\hat{p}(\hat{x})) \\
&=\sum_{x\in p\Fix(\beta\alpha\tilde{f}_i)}|p'\Fix(\beta\alpha\tilde{f}_i)\cap \hat{p}^{-1}(x)|\ind(f;x) \\
&=\sum_{x\in p\Fix(\beta\alpha\tilde{f}_i)}|u_i(\coin(\tau_{\beta\alpha}\varphi_i,\iota_i))|\ind(f;x) \\
&=|u_i(\coin(\tau_{\beta\alpha}\varphi_i,\iota_i))|\ind(f;p\Fix(\beta\alpha\tilde{f}_i)). \qedhere
\end{align*}
\end{proof}

With this, our proof of Theorem \ref{thm:av-RT} is finished.

\section{Averaging formula for the Lefschetz number}

Using the results from the previous section, we immediately obtain:

\begin{thm}\label{thm:av-L}
For an $n$-valued map $f:X\to D_n(X)$ with lift $(\tilde{f}_1,\ldots,\tilde{f}_n):\tilde{X}\to F_n(\tilde{X},\pi)$ and corresponding split lift $(\bar{f}_1,\ldots,\bar{f}_n)$ with respect to $S$ and $\Gamma$, \[
L(f)=\frac{1}{[\pi:S]}\sum_{i=1}^{n}\sum_{\bar{\alpha}\in \pi/\Gamma} L(\bar{\alpha}\bar{f}_i,q).
\]
\end{thm}

\begin{proof}
We have \[
L(f) = \sum_{[(\alpha,i)]\in \R[\tilde{f}_\#]} \ind(f;p\Fix(\alpha\tilde{f}_i)).
\]
By applying Lemma \ref{lem:decomp} to the terms \[
X_{[(\alpha,i)]}=\ind(f;p\Fix(\alpha\tilde{f}_i))
\]
in the abelian group $\ZZ$, this sum can be decomposed as \[
L(f) = \frac{1}{[\pi:S]}\sum_{i=1}^{n}\sum_{\bar{\alpha}\in \pi/\Gamma}\sum_{[\beta]\in \R[\tau_\alpha\varphi'_i,\iota]} |u_i(\coin(\tau_{\beta\alpha}\varphi_i,\iota_i))|\ind(f;p\Fix(\beta\alpha\tilde{f}_i)).
\] 
By Lemma \ref{lem:ind}, this becomes \begin{align*}
L(f) &= \frac{1}{[\pi:S]}\sum_{i=1}^{n}\sum_{\bar{\alpha}\in \pi/\Gamma}\sum_{[\beta]\in \R[\tau_\alpha\varphi'_i,\iota]} \ind(\bar{\alpha}\bar{f}_i,q;p'\!\Coin(\beta\alpha\tilde{f}_i,\id)) \\
&= \frac{1}{[\pi:S]}\sum_{i=1}^{n}\sum_{\bar{\alpha}\in \pi/\Gamma} L(\bar{\alpha}\bar{f}_i,q). \qedhere
\end{align*}
\end{proof}

\section{Averaging formula for the Nielsen number}

As in the single-valued case, only a partial averaging formula can be obtained for the Nielsen number. Analogous to the main result of \cite{kimleelee}, we get:

\begin{thm}\label{thm:av-N}
For an $n$-valued map $f:X\to D_n(X)$ with lift $(\tilde{f}_1,\ldots,\tilde{f}_n):\tilde{X}\to F_n(\tilde{X},\pi)$ and corresponding split lift $(\bar{f}_1,\ldots,\bar{f}_n)$ with respect to $S$ and $\Gamma$, \[
N(f)\geq \frac{1}{[\pi:S]}\sum_{i=1}^{n}\sum_{\bar{\alpha}\in \pi/\Gamma} N(\bar{\alpha}\bar{f}_i,q)
\]
with an equality if and only if $\coin(\tau_{\beta\alpha}\varphi_i,\iota_i)\subseteq S$ whenever $p'\!\Coin(\beta\alpha \tilde{f}_i,\id)$ is an essential coincidence class of $(\bar{\alpha}\bar{f}_i,q)$.
\end{thm}

\begin{proof}
For a fixed point class $\FF$ of a map $f$, write \[
\eps(f;\FF)=\left\{
\begin{array}{ll}
1 & \text{if $\ind(f;\FF)\neq 0$}, \\
0 & \text{otherwise};
\end{array}
\right.
\]
and similarly for coincidences. Then we can express \[
N(f)=\sum_{[(\alpha,i)]\in \R[\tilde{f}_\#]} \eps(f;p\Fix(\alpha\tilde{f}_i))
\]
and \[
N(\bar{\alpha}\bar{f}_i,q)=\sum_{[\beta]\in \R[\tau_\alpha\varphi'_i,\iota]} \eps(\bar{\alpha}\bar{f}_i,q;p'\!\Coin(\beta\alpha\tilde{f}_i,\id)).
\]
By applying Lemma \ref{lem:decomp} to the terms \[
X_{[(\alpha,i)]}=\eps(f;p\Fix(\alpha\tilde{f}_i))
\]
in the abelian group $\ZZ$, we can decompose \[
N(f) = \frac{1}{[\pi:S]}\sum_{i=1}^{n}\sum_{\bar{\alpha}\in \pi/\Gamma}\sum_{[\beta]\in \R[\tau_\alpha\varphi'_i,\iota]} |u_i(\coin(\tau_{\beta\alpha}\varphi_i,\iota_i))|\,\eps(f;p\Fix(\beta\alpha\tilde{f}_i)).
\]
By Lemma \ref{lem:ind}, we have \[
\eps(f;p\Fix(\beta\alpha\tilde{f}_i))=\eps(\bar{\alpha}\bar{f}_i,q;p'\!\Coin(\beta\alpha\tilde{f}_i,\id))
\]
for all $\alpha$, $\beta$ and $i$, so that \begin{align*}
N(f) &= \frac{1}{[\pi:S]}\sum_{i=1}^{n}\sum_{\bar{\alpha}\in \pi/\Gamma}\sum_{[\beta]\in \R[\tau_\alpha\varphi'_i,\iota]} |u_i(\coin(\tau_{\beta\alpha}\varphi_i,\iota_i))|\,\eps(\bar{\alpha}\bar{f}_i,q;p'\!\Coin(\beta\alpha\tilde{f}_i,\id)) \\
&\geq \frac{1}{[\pi:S]}\sum_{i=1}^{n}\sum_{\bar{\alpha}\in \pi/\Gamma}\sum_{[\beta]\in \R[\tau_\alpha\varphi'_i,\iota]} \eps(\bar{\alpha}\bar{f}_i,q;p'\!\Coin(\beta\alpha\tilde{f}_i,\id)) \\
&= \frac{1}{[\pi:S]}\sum_{i=1}^{n}\sum_{\bar{\alpha}\in \pi/\Gamma} N(\bar{\alpha}\bar{f}_i,q),
\end{align*}
with an equality if and only if $u_i(\coin(\tau_{\beta\alpha}\varphi_i,\iota_i))=1$, i.e.\ $\coin(\tau_{\beta\alpha}\varphi_i,\iota_i)\subseteq S$, whenever $\eps(\bar{\alpha}\bar{f}_i,q;p'\!\Coin(\beta\alpha\tilde{f}_i,\id))\neq 0$.
\end{proof}

\section{Infra-nilmanifolds}\label{sec:inm}

Let $G$ be a connected and simply connected nilpotent Lie group. A nilmanifold is a quotient $\Gamma\orb G$, where $\Gamma\subseteq G$ is a lattice (i.e.\ a discrete and cocompact subgroup). An infra-nilmanifold is a compact manifold that is finitely covered by a nilmanifold. Equivalently, an infra-nilmanifold can be expressed as a quotient $\pi\orb G$, where $\pi\subseteq G\rtimes \Aut(G)$ is a so-called almost-Bieberbach group: a torsion free cocompact discrete subgroup of $G\rtimes C$, with $C$ a maximal compact subgroup of $\Aut(G)$. In that case, an example of a nilmanifold that finitely covers $\pi\orb G$ is $\Gamma\orb G$, where $\Gamma=\pi\cap G$ is the subgroup of `pure translations'.

Lefschetz and Nielsen coincidence numbers on nilmanifolds have been studied in \cite{HLP2011,HLP2012}: for two maps $f,g:\Gamma\orb G\to \Gamma'\orb G'$ between equal-dimensional nilmanifolds, let $\varphi,\psi:\Gamma\to \Gamma'$ be their induced morphisms. These morphisms extend uniquely to Lie group morphisms $G\to G'$; if $\varphi_*,\psi_*:\g\to\g'$ denote the corresponding morphisms on the Lie algebras of $G$ and $G'$, \[
L(f,g)=\det_{\Gamma}^{\Gamma'}(\psi_*-\varphi_*) \qquad \text{and} \qquad
N(f,g)=\left|\det_{\Gamma}^{\Gamma'}(\psi_*-\varphi_*)\right|
\]
were $\displaystyle\det_{\Gamma}^{\Gamma'}$ denotes the determinant of the matrix
with respect to certain specific bases of $\g$ and $\g'$ induced by $\Gamma$ and $\Gamma'$; see \cite{HLP2012}. We will only need the special case where $G=G'$ and $\Gamma$ is a finite index subgroup of $\Gamma'$, in which \cite[Corollary 6.6]{HLP2012} can be used to simplify \[
\det_{\Gamma}^{\Gamma'}(\psi_*-\varphi_*)=[\Gamma':\Gamma]\det_{\Gamma}^{\Gamma}(\psi_*-\varphi_*)=[\Gamma':\Gamma]\det(\psi_*-\varphi_*).
\]
Here the last determinant is the usual determinant for morphisms $\g\to\g$, defined unambiguously by taking the same (arbitrary) basis for the domain and image.

Together with our results, we get:

\begin{thm}
Suppose $f:\pi\orb G\to D_n(\pi\orb G)$ is an $n$-valued map on an infra-nilmanifold $\pi\orb G$, finitely covered by a nilmanifold $\Gamma\orb G$. Let $S$ be an $f$-$\Gamma$-invariant subgroup of $\pi$. Let $\varphi_i:S_i\to \pi$ be the morphisms induced by some lift $\tilde{f}$ of $f$, with restrictions $\varphi'_i:S\to \Gamma$. For all $\alpha$ and $i$, let $(\tau_\alpha\varphi'_i)_*:\g\to\g$ denote the Lie algebra morphism induced by the unique extension of $\tau_\alpha\varphi'_i:S\to \Gamma$ to a morphism $G\to G$. Then \[
L(f)=\frac{1}{[\pi:\Gamma]}\sum_{i=1}^{n}\sum_{\bar{\alpha}\in \pi/\Gamma} \det(I-(\tau_\alpha\varphi'_i)_*),
\]
where $I:\g\to\g$ is the identity Lie algebra morphism.
\end{thm}

\begin{proof}
Let $(\bar{f}_1,\ldots,\bar{f}_n)$ be a split lift of $f$ with respect to $S$ and $\Gamma$. As we saw at the very end of section \ref{sec:split-lift}, the morphism induced by $\bar{\alpha}\bar{f}_i$ is $\tau_\alpha\varphi'_i:S\to \Gamma$. Since $S$ is a finite index subgroup of $\Gamma$, we can apply the above result to write \[
L(\bar{\alpha}\bar{f}_i,q)=[\Gamma:S]\det(\iota_*-(\tau_\alpha\varphi'_i)_*).
\]
Since the unique extension $G\to G$ of the inclusion $\iota:S\to \Gamma$ is the identity, whose induced Lie algebra morphism is the identity $I:\g\to\g$, this becomes \[
L(\bar{\alpha}\bar{f}_i,q)=[\Gamma:S]\det(I-(\tau_\alpha\varphi'_i)_*).
\]
Plugging this into Theorem \ref{thm:av-L} gives \begin{align*}
L(f) &= \frac{1}{[\pi:S]}\sum_{i=1}^{n}\sum_{\bar{\alpha}\in \pi/\Gamma} [\Gamma:S]\det(I-(\tau_\alpha\varphi'_i)_*) \\
&= \frac{1}{[\pi:\Gamma]}\sum_{i=1}^{n}\sum_{\bar{\alpha}\in \pi/\Gamma} \det(I-(\tau_\alpha\varphi'_i)_*). \qedhere
\end{align*}
\end{proof}

Note that in this setting, an $f$-$\Gamma$-invariant subgroup always exists, by Remark \ref{rmk:virt-polyc} (indeed, finite extensions of finitely generated nilpotent groups are virtually polycyclic).

\medskip

As in the single-valued case, the inequality for the Nielsen number is always an equality for infra-nilmanifolds:

\begin{lem}
In the setting of the previous sections, if $X=\pi\orb G$ is an infra-nilmanifold and $\bar{X}=\Gamma\orb G$ a nilmanifold, if $\ind(\bar{\alpha}\bar{f}_i,q;p'\!\Coin(\beta\alpha\tilde{f}_i,\id))\neq 0$, then $\coin(\tau_{\beta\alpha}\varphi_i,\iota_i)=1$.
\end{lem}

\begin{proof}
If the index $\ind(\bar{\alpha}\bar{f}_i,q;p'\!\Coin(\beta\alpha\tilde{f}_i,\id))$ is non-zero, then so is the Nielsen number \[
N(\bar{\alpha}\bar{f}_i,q)=N(\overline{\beta\alpha}\bar{f}_i,q)=[\Gamma:S]|\det(I-(\tau_{\beta\alpha}\varphi'_i)_*)|.
\]
In particular, in that case, $\det(I-(\tau_{\beta\alpha}\varphi'_i)_*)\neq 0$. It follows from \cite[Lemma 4.1]{affien}\footnote{To be precise, not from this lemma as stated in \cite{affien}, but from its more general version where $S'$ is allowed to be any finite index subgroup; the proof remains perfectly valid in that case.} that $\coin(\tau_{\beta\alpha}\varphi_i,\iota_i)=1$.
\end{proof}

Similarly as for $L(f)$, we get:

\begin{thm}\label{thm:av-N-inm}
Suppose $f:\pi\orb G\to D_n(\pi\orb G)$ is an $n$-valued map on an infra-nilmanifold $\pi\orb G$, finitely covered by a nilmanifold $\Gamma\orb G$. Let $S$ be an $f$-$\Gamma$-invariant subgroup of $\pi$. Let $\varphi_i:S_i\to \pi$ be the morphisms induced by some lift $\tilde{f}$ of $f$, with restrictions $\varphi'_i:S\to \Gamma$. For all $\alpha$ and $i$, let $(\tau_\alpha\varphi'_i)_*:\g\to\g$ denote the Lie algebra morphism induced by the unique extension of $\tau_\alpha\varphi'_i:S\to \Gamma$ to a morphism $G\to G$. Then \[
N(f)=\frac{1}{[\pi:\Gamma]}\sum_{i=1}^{n}\sum_{\bar{\alpha}\in \pi/\Gamma} |\det(I-(\tau_\alpha\varphi'_i)_*)|.
\]
\end{thm}

\begin{rmk}
Three choices are involved in computing the right hand sides of our formulas for $L(f)$ and $N(f)$:
\begin{itemize}
\item[(i)] The choice of $f$-$\Gamma$-invariant subgroup $S$;
\item[(ii)] The choice of reference lift $\tilde{f}$ for $f$, which determines the morphisms $\varphi_i$;
\item[(iii)] The choice of reference lift $\alpha\tilde{f}_j$ for $\bar{\alpha}\bar{f}_j$, i.e.\ the choice of representative $\alpha\in \pi$ for $\bar{\alpha}\in \pi/\Gamma$, which appears in $\det(I-(\tau_\alpha\varphi'_i)_*)$.
\end{itemize}
As a sanity check, let us show that the final formulas are indeed independent of these choices.

For (i), let $S$ and $S'$ be two $f$-$\Gamma$-invariant subgroups of $\pi$, and $\varphi'_i:S\to \Gamma$ and $\varphi''_i:S'\to \Gamma$ the corresponding morphisms. Since both are restrictions of the same morphism $\varphi_i:S_i\to \pi$, they are the same on $S\cap S'$, which is also a lattice in $G$. Both the extensions of $\varphi'_i$ and $\varphi''_i$ to Lie group morphisms $G\to G$ are also extensions of this morphism $S\cap S'\to \Gamma$. By uniqueness, they must coincide.

For (ii), recall from Remark \ref{rmk:lift-choice} that, if $\varphi_i:S_i\to \pi$ are the morphisms induced by one lift of $f$, another lift of $f$ induces morphisms $\tau_{\gamma_i}\varphi_{\eta^{-1}(i)}:S_{\eta^{-1}(i)}\to \pi$ for some $\gamma_i\in \pi$ and $\eta\in\S_n$. 
For these fixed $\gamma_i\in \pi$ and $\eta\in\S_n$, as $\bar{\alpha}$ runs over all elements of $\pi/\Gamma$, so does $\bar{\beta}=\overline{\alpha\gamma_i}$; and as $i$ runs over $1,\ldots,n$, so does $j=\eta^{-1}(i)$. Thus, the formula we get for $L(f)$ becomes \begin{align*}
\frac{1}{[\pi:\Gamma]} \sum_{i=1}^n \sum_{\bar{\alpha}\in \pi/\Gamma} \det(I-(\tau_{\alpha\gamma_i}\varphi'_{\eta^{-1}(i)})_*) &= \frac{1}{[\pi:\Gamma]} \sum_{i=1}^n \sum_{\bar{\beta}\in \pi/\Gamma} \det(I-(\tau_\beta\varphi'_{\eta^{-1}(i)})_*) \\
&= \frac{1}{[\pi:\Gamma]} \sum_{j=1}^n \sum_{\bar{\beta}\in \pi/\Gamma} \det(I-(\tau_\beta\varphi'_j)_*).
\end{align*}
Analogously, the formula for $N(f)$ is independent of the chosen lift $\tilde{f}$ as well.

For (iii), we recall the following result.

\begin{lem}[{\cite[Lemma 3.2]{leelee2006}}]\label{lem:leelee}
Let $G$ be a connected simply connected nilpotent Lie group with Lie algebra $\g$. For any morphism $\varphi:G\to G$ and for any $g\in G$, \[
\det(I-(\tau_g)_*\varphi_*)=\det(I-\varphi_*).
\]
\end{lem}

If $\beta\alpha$ is another representative for $\bar{\alpha}\in \pi/\Gamma$, with $\beta\in \Gamma\subseteq G$, we can take $\varphi$ in this lemma to be the unique extension of $\tau_{\alpha}\varphi'_i:S\to \Gamma$ and $g=\beta$, and we get \[
\det(I-(\tau_{\beta\alpha}\varphi'_i)_*)=\det(I-(\tau_\beta)_*(\tau_{\alpha}\varphi'_i)_*)=\det(I-(\tau_\alpha\varphi'_i)_*).
\]
\end{rmk}

\begin{rmk}
Also using Lemma \ref{lem:leelee}, the formula for the Nielsen number of an affine $n$-valued map on an infra-nilmanifold from \cite{affien} can be obtained as a special case of Theorem \ref{thm:av-N-inm}.
\end{rmk}

\begin{ex}
To illustrate our formula, we can apply it to the map $f$ from Example \ref{ex:klein}.

Recall that we identified the Klein bottle with $\pi\orb \RR^2$, where $\pi\subseteq \RR^2\rtimes \GL_2(\RR)$ is the group generated by \[
a=\left(\begin{bmatrix}1 \\ 0 \end{bmatrix},\begin{bmatrix} 1 & 0 \\ 0 & 1 \end{bmatrix}\right), \quad b=\left(\begin{bmatrix}0 \\ \frac{1}{2} \end{bmatrix},\begin{bmatrix} -1 & 0 \\ 0 & 1 \end{bmatrix}\right).
\]
As a finite cover for $\pi\orb \RR^2$, we can take the torus $\Gamma\orb \RR^2$, with $\Gamma=\langle a,b^2\rangle \cong \ZZ^2$. Recall that the morphism $\tilde{f}_\#:\pi\to \pi^2\rtimes \S_2$ induced by the lift $\tilde{f}$ from Example \ref{ex:klein} is given by $\tilde{f}_\#(a)=(1,1;1)$ and $\tilde{f}_\#(b)=(b,1;\sigma)$ where $\sigma$ is the permutation interchanging $1$ and $2$, so $S_1=S_2=\langle a,b^2\rangle$ and the morphisms $\varphi_i:S_i\to \pi$ are \[
\varphi_1=\varphi_2: \langle a,b^2\rangle\to \pi:a^kb^{2\ell}\mapsto b^\ell.
\]
As an $f$-$\Gamma$-invariant subgroup we can take \[
S=\langle a,b^4\rangle\subseteq \pi.
\]
The corresponding morphisms $\varphi'_i:S\to \Gamma$ are \[
\varphi'_1=\varphi'_2: S\to \Gamma:a^k b^{4\ell}\mapsto b^{2\ell}.
\]
As representatives in $\pi$ for the two elements of $\pi/\Gamma$, we can take $1$ and $b$, with \begin{align*}
\tau_1\varphi'_i(a^kb^{4\ell})&=b^{2\ell} \\
\tau_b\varphi'_i(a^kb^{4\ell})&=bb^{2\ell}b^{-1}=b^{2\ell},
\end{align*}
or, under the isomorphism $\Gamma\cong \ZZ^2$, \[
\tau_1\varphi'_i\left(\begin{bmatrix} k \\ 2\ell \end{bmatrix}\right)=\tau_b\varphi'_i\left(\begin{bmatrix} k \\ 2\ell \end{bmatrix}\right)=\begin{bmatrix} 0 \\ \ell \end{bmatrix}.
\]
The unique extensions $\RR^2\to \RR^2$ (and also their Lie algebra morphisms) are given by the matrix \[
(\tau_1\varphi'_i)_*=\begin{bmatrix} 0 & 0 \\ 0 & \frac{1}{2} \end{bmatrix}=(\tau_b\varphi'_i)_*.
\]
We get \begin{align*}
L(f) &= \frac{1}{2} \sum_{i=1}^{2} \sum_{\bar{\alpha}\in \pi/\Gamma} \det\left(\begin{bmatrix} 1 & 0 \\ 0 & 1 \end{bmatrix}-\begin{bmatrix} 0 & 0 \\ 0 & \frac{1}{2} \end{bmatrix}\right) = \frac{1}{2}\cdot 2\cdot 2\cdot \frac{1}{2}=1
\end{align*}
and similarly
\begin{align*}
N(f) &= \frac{1}{2} \sum_{i=1}^{2} \sum_{\bar{\alpha}\in \pi/\Gamma} \left|\det\left(\begin{bmatrix} 1 & 0 \\ 0 & 1 \end{bmatrix}-\begin{bmatrix} 0 & 0 \\ 0 & \frac{1}{2} \end{bmatrix}\right)\right| =1,
\end{align*}
which is consistent with our result from \cite[Example 5.1]{affien}.
\end{ex}

\end{document}